\documentclass[final,twoside,11pt]{entics} 
\usepackage{enticsmacro}

\usepackage{tikz}

\usepackage{graphicx}
\usepackage[all]{xy}
\newcommand{\bjoin}{\bigvee}
\newcommand{\join}{\vee}
\newcommand{\bmeet}{\bigwedge}
\newcommand{\meet}{\wedge}

\newcommand{\ctin}{\subseteq}

\newcommand{\da}{\mbox{$\downarrow$}}

\newcommand{\bproof}{\vskip-10pt\begin{pf}\hskip-3pt
	}
	\newcommand{\eproof}{\end{pf}}

\newcommand{\calA}{\mathcal{A}}

\newcommand{\calS}{\mathcal{S}}
\newcommand{\imp}{\Rightarrow}

\newcommand{\dimp}{\Longleftrightarrow}
\newcommand{\limp}{\Leftarrow}

\newcommand{\benu}{\begin{enumerate}}
	\newcommand{\eenu}{\end{enumerate}}
\newcommand{\blfl}{\begin{list}{}}
	\newcommand{\elfl}{\end{list}}

\newcommand{\np}{\vskip5pt\noindent}

\newcommand{\darw}{\operatorname{\downarrow}}

\newcommand{\rnum}{\mathbb{R}}
\newcommand{\qnum}{\mathbb{Q}}

\newcommand{\sels}{\mathcal{S}}

\newcommand{\sframe}{$\calS$-frame}

\newcommand{\aichs}{\mathcal{H}_{\mathcal{S}}}

\newcommand{\sideal}{$\mathcal{S}$-ideal}
\newcommand{\sfrm}{$\sels$\textbf{Frm}}

\newcommand{\bjoinl}[1]{\bjoin\limits_{#1}}
\newcommand{\bcupl}[1]{\bigcup\limits_{#1}}

\newcommand{\slind}{$\sels$-Lindel\"{o}f}

\newcommand{\sframes}{\sframe{}s}
\newcommand{\myopdar}{\downarrow\!}

\newcommand{\calC}{\mathcal{C}}
\newcommand{\calP}{\mathcal{P}}

\newcommand{\scong}{$\sels$-congruence}
\newcommand{\congs}{\mathcal{C}_{\sels}}
\newcommand{\congsL}{\mathcal{C}_{\sels}L}
\newcommand{\congsM}{\mathcal{C}_{\sels}M}

\newcommand{\naba}{\nabla\hskip-2pt_a}
\newcommand{\dela}{\Delta_a}

\newcommand{\aichsL}{\mathcal{H}_{\mathcal{S}}L}

\newcommand{\sym}{\raisebox{\depth}{$\chi$}}

\newcommand{\chio}{\sym_0}

\newcommand{\msl}{meet-semilattice}

\newcommand{\bnote}{\begin{note}}
	\newcommand{\enote}{\end{note}}
\newcommand{\brem}{\begin{remark}}
	\newcommand{\erem}{\end{remark}}
\newcommand{\bdefn}{\begin{definition} }
	\newcommand{\edefn}{\end{definition}}
\newcommand{\blemma}{\begin{lemma}}
	\newcommand{\elemma}{\end{lemma}}
\newcommand{\bthm}{\begin{theorem}}
	\newcommand{\ethm}{\end{theorem}}
\newcommand{\bcrol}{\begin{corollary}}
	\newcommand{\ecrol}{\end{corollary}}
\newcommand{\bexa}{\begin{example}}
	\newcommand{\eexa}{\hfill\qed\end{example}}
\newcommand{\bprop}{\begin{proposition}}
	\newcommand{\eprop}{\end{proposition}}

\def\conf{ISDT 2022} 	
\volume{2}			
\def\volu{2}			
\def\papno{10}			
\def\mydoi{{\normalshape  \href{https://doi.org/10.46298/entics.10459}{10.46298/entics.10459}}}
\def\copynames{A. Schauerte, J. Frith} 
\def\runtitle{Partial frames, their free frames and their congruence frames}

\def\lastname{Schauerte snd Firth}
\begin{document}
\begin{frontmatter}
  \title{Partial Frames, Their Free Frames\\ and Their Congruence Frames} 
  \author{Anneliese Schauerte\thanksref{myemail}}	
   \author{John Frith\thanksref{coemail}}		
  \thanks[myemail]{Email: \href{mailto:Anneliese.Schauerte@uct.ac.za} {\texttt{\normalshape     Anneliese.Schauerte@uct.ac.za}}} 
  \address[b]{Department of Mathematics and Applied Mathematics\\University of Cape Town\\
    Cape Town, South Africa} 
  \thanks[coemail]{Email:  \href{mailto:John.Frith@uct.ac.za} {\texttt{\normalshape
        John.Frith@uct.ac.za}}}
\begin{abstract} 
 The context of this work is that of partial frames; these  are meet-semilattices where not all subsets need have joins. A selection function, $\sels$, specifies, for all meet-semilattices, certain subsets under consideration, which we call the ``designated'' ones;  an \sframe{} then must have joins of (at least) all such subsets and binary meet must distribute over these.  A small collection of axioms suffices to specify our selection functions; these axioms are sufficiently general to  include as examples of partial frames, bounded distributive lattices, $\sigma$-frames,  $\kappa$-frames and frames.
 
 We consider  right and left adjoints of \sframe{} maps, as a prelude to the introduction of closed and open maps. 
 
 Then we look at what might be an appropriate notion of Booleanness for partial frames. The obvious candidate is the condition that  every element be complemented; this concept is indeed of interest, but we pose three further conditions which, in the frame setting, are all equivalent to it. However, in the context of partial frames, the four conditions are distinct. In investigating these, we make essential use of the free frame over a partial frame and the congruence frame of a partial frame.
 
 We  compare congruences of a partial frame, technically called \scong{s}, with the frame congruences of its free frame. We provide a natural transformation for the situation and also consider right adjoints of the frame maps in question. We characterize the case where the two congruence frames are isomorphic and provide examples which illuminate the possible different behaviour of the two. 
 
 We conclude with a characterization of closedness and openness for the embedding of a partial frame into its free fame, and into its congruence frame. 
\end{abstract}
\end{frontmatter}
\section{Introduction}\label{intro}

Partial frames are meet-semilattices where, in contrast with frames, not all subsets need have joins. A selection function, $\sels$, specifies, for all meet-semilattices, certain subsets under consideration, which we call the ``designated'' ones;  an \sframe{} then must have joins of (at least) all such subsets and binary meet must distribute over these.  A small collection of axioms suffices to specify our selection functions; these axioms are sufficiently general to
include as examples of partial frames, bounded distributive lattices, $\sigma$-frames,
$\kappa$-frames and frames.

We consider the classical notions of right and left adjoints, for \sframe{} maps. Unlike the situation for full frames, such maps need not have right adjoints. This is a prelude to the introduction of closed and open maps, and a discussion of their properties. 

What is an appropriate notion of Booleanness for partial frames? The obvious answer is that the partial frame should have every element complemented; this concept is indeed of interest, but we pose three further conditions which, in the frame setting, are all equivalent to it. However, in the context of partial frames, the four conditions are distinct. In investigating these, we make essential use of the free frame over a partial frame and the congruence frame of a partial frame.

We  compare congruences of a partial frame, technically called \scong{s}, with the frame congruences of its free frame. We provide a natural transformation for the situation and also consider right adjoints of the frame maps in question. We characterize the case where the two congruence frames are isomorphic and provide examples which illuminate the possible different behaviour of the two. 

We conclude with a characterization of closedness and openness for the embedding of a partial frame into its free fame, and into its congruence frame. 

Since this document is intended as an extended abstract, proofs are omitted.
\section{Background}

This background section is taken largely from \cite{jfasBoolVar}. 
See \cite{picpulframesbook} and \cite{jo}  as
references for frame theory; see \cite{BBCrag} and \cite{bbnotes} for $\sigma$-frames; see \cite{madden} and \cite{manpap} for $\kappa$-frames; see \cite{maclane} and \cite{ahs} for general category theory.

\np The basics of our approach to partial frames can be found in \cite{jfascgasa1}, \cite{jfascgasa2} and \cite{jfas_hcozforqm}. Our papers with a more topological flavour are \cite{jfascplns},  \cite{jfasslovak},  \cite{jfasonepoint}, \cite{jfaspartframfilt} and \cite{jfascompreflpart}.  Our papers with a more algebraic flavour are \cite{jfascov}, \cite{jfasmadden} and \cite{jfassemilatticescong}. Crucial for this paper is  \cite{jfasheyt}. We are indebted to earlier work by other authors in this field: see \cite{paseka},  \cite{zhaonuclei}, \cite{zhao} and \cite{zenk}. For those interested in a comparison of the various approaches, see \cite{jfascgasa2}.
\np
A \textit{meet-semilattice} is a partially ordered set in which all finite subsets have a meet. In particular,
we regard the empty set as finite, so a meet-semilattice comes equipped with a top element, which we denote by $1$.
We do not insist that a meet-semilattice should have a bottom element, which, if it exists, we denote by $0$.
A function between meet-semilattices $f:L\to M$  is a \textit{meet-semilattice map} if it preserves finite meets, as well as the top element. A \textit{sub \msl{}} is a subset for which the inclusion map is a \msl{} map.

The essential idea for a \textit{partial frame} is that it should be ``frame-like'' but that not all joins need exist; only certain joins have guaranteed existence and binary meets should distribute over these joins. The guaranteed joins are specified in a global way on the category of \msl{s} by specifying what is called a selection function; the details are given below.
\begin{definition}\label{AB}
A \textit{selection function} is a rule, which we usually denote by $\sels$, which assigns to each meet-semilattice $A$ a collection $\sels A$ of subsets of $A$, such that the following conditions hold (for all meet-semilattices $A$ and $B$):
\begin{enumerate}
\item[] 
\begin{enumerate}
\item[(S1)] For all $x\in A$, $\{x\}\in\sels A$.
\item[(S2)] If $G, H\in\sels A$ then $\{x\meet y:x\in G,y\in H\}\in\sels A$.
\item[(S2)$'$] If $G,H\in\sels A$ then $\{x\join y:x\in G,y\in H\}\in\sels A$.
\item[(S3)] If $G\in\sels A$ and, for all $x\in G$, $x=\bjoin H_x$ for some $H_x\in \sels A$, then $$\bigcup\limits_{x\in G}H_x\in\sels A.$$
\item[(S4)] For any meet-semilattice map $f:A\to B$,
$$\sels(f[A])=\{f[G]:G\in\sels A\}\ctin \sels B.$$
\item[(SSub)] For any sub \msl{} $B$ of \msl{} $A$, if $G\ctin B$ and $G\in\sels A $, then $G\in\sels B$.
\item[(SFin)] If $F$ is a finite subset of $A$, then $F\in\sels A$. 
\item[(SCov)] If $G\ctin H$ and $H\in\sels A$ with $\bjoin H=1$ then $G\in\sels A$. (Such $H$ are called $\sels$-\textit{covers}.)
\item[(SRef)] Let $X,Y\ctin A$. If $X\leq Y$ with $X\in\sels A$  there is a $C\in \sels A$ such that $X\leq C\ctin Y$. (By $X\leq Y$ we mean, as usual, that for each $x\in X$ there exists $y\in Y$ such that $x\leq y$.)
\end{enumerate}
\end{enumerate}
\end{definition}
Of course (SFin) implies (S1) but there are situations where we do not impose (SFin) but insist on (S1). Note that we always have $\emptyset\in\sels A$.
Once a selection function, $\sels$, has been fixed, we speak informally of the members of $\sels A$ as the \textit{designated} subsets of $A$.

\begin{definition}
An \textit{\sframe}  $L$ is a \msl{} in which every designated subset has a join and for any such designated subset $B$ of $L$ and any $a\in L$, $$a\meet\bjoin B=\bjoinl{b\in B}a\meet b.$$Of course such an \sframe{} has both a top and a bottom element which we denote by $1$ and $0$ respectively.\\
A \msl{} map $f:L\to M$, where $L$ and $M$ are \sframe{s}, is an \textit{\sframe{} map} if $f(\bjoin B)=\bjoinl{b\in B}f(b)$ for any designated subset $B$ of $L$. In particular such an $f$ preserves the top and bottom element.\\
A \textit{sub \sframe{}} $T$ of an \sframe{} $L$ is a subset of $L$ such that the inclusion map $i:T\to L$ is an \sframe{} map. \\
The category \sfrm{} has objects \sframe{s} and arrows \sframe{} maps.
\end{definition}

We use the terms ``partial frame'' and ``\sframe'' interchangeably, especially if no confusion about the selection function is likely. We also use the term \textit{full frame} in situations where we wish to emphasize that all joins exist.

\begin{note} Here are some examples of different selection functions and their corresponding \sframe{s}.
\begin{itemize}
\item[1.] In the case that all joins are specified, we are of course considering the notion of a frame. 
\item[2.] In the case that (at most) countable joins are specified, we have the notion of a $\sigma$-frame. 
\item[3.] In the case that joins of subsets with cardinality less than some (regular) cardinal $\kappa$ are specified, we have the notion of a $\kappa$-frame.
\item[4.] In the case that only finite joins are specified, we have the notion of a bounded distributive lattice.
\end{itemize}
\end{note}

The remainder of this section gives a lot of information about $\aichsL$, the free frame over the \sframe{} $L$, as well as $\congsL$, the frame of \scong{s} of $L$, and the relationship between the two. These results come from \cite{jfas_hcozforqm} on $\aichsL$, \cite {jfascov} and \cite{jfasmadden} on $\congsL$.

In the definition below, $L$ is an \sframe{}. 

\begin{definition}\label{aaaa}
\begin{enumerate}\item[(a)] A subset $J$ of an  $L$ is an \textit{$\mathcal{S}$-ideal of $L$}
if $J$ is a non-empty downset closed under  designated joins (the latter
meaning that if $X\ctin J$, for $X$ a designated subset of $L$, then $\bjoin
X\in J$).
\item[(b)] The collection of all \sideal s of  $L$ will be denoted
$\aichs L$, and called the \textit{\sideal{} frame of $L$}. It is in fact the free frame over $L$. 
\item[(c)] For $I\in \aichsL$, $t\in(\darw x)\join I\dimp t\leq x\join s$, for some $s\in I$.
\item[(d)] We call $\theta\ctin L\times L$ an \textit{\scong{} on }$L$ if it satisfies the following:\\
(C1) $\theta$ is an equivalence relation on $L$.\\
(C2) $(a,b), (c,d)\in\theta$ implies that $(a\meet c,b\meet d)\in\theta$.\\
(C3) If $\{(a_{\alpha},b_{\alpha}):\alpha\in\calA\}\ctin\theta$ and $\{a_{\alpha}:\alpha\in\calA\}$ and $\{b_{\alpha}:\alpha\in\calA\}$ are designated subsets of $L$, then $(\bjoinl{\alpha\in\calA}a_{\alpha},\bjoinl{\alpha\in\calA}b_{\alpha})\in\theta$.
\item[(e)] The collection of all \scong{s} on $L$ is denoted by $\congsL$; it is in fact a (full) frame with meet given by intersection. 
\item[(f)] 
\begin{enumerate}
	\item[(i)] Let $A\ctin L\times L$. We use the notation $\langle A\rangle$ to denote the smallest \scong{} containing $A$. This exists by completeness of $\congsL$.
	\item[(ii)] We define $\naba=\{(x,y):x\join a=y\join a\}$ and  $\Delta_a=\{(x,y):x\meet a=y\meet a\}$; these are \scong{s} on $L$.
	\item[(iii)] It is easily seen that $\naba=\bigcap\{\theta:\theta\in\congsL\textrm{ and }(0,a)\in\theta\}=\langle(0,a)\rangle$ and that $\Delta_a=\bigcap\{\theta:\theta\in\congsL\textrm{ and }(a,1)\in\theta\}=\langle(a,1)\rangle$. 
	\item[(iv)] For $a\leq b$, it follows that $\Delta_a\cap\nabla_b=\langle(a,b)\rangle$. 
	\item[(v)] The congruence $\nabla_1=L\times L$ is the top element and $\nabla_0=\{(x,x):x\in L\}$ (called the \textit{diagonal}) is the bottom element of $\congsL$.
\end{enumerate}

\item[(g)] The following hold in $\congsL$.
\benu
\item[(i)] For any $\theta\in\congsL$, $\theta=\bjoin\{\nabla_b\meet\Delta_a:(a,b)\in\theta,a\leq b\}$.
\item[(ii)] $\naba\join\theta=\{(x,y):(x\join a,y\join a)\in\theta\}$.
\item[(iii)] $\dela\join\theta=\{(x,y):(x\meet a,y\meet a)\in\theta\}$.
\item[(iv)] For any $I\in\aichsL$, $\bjoinl{x\in I}\nabla_x=\bcupl{x\in I}\nabla_x$.
\eenu
\item[(h)] The function $\nabla:L\to\congsL$ given by $\nabla(a)=\nabla_a$ is an \sframe{} embedding. It has the universal property that if $f:L\to M$ is an \sframe{} map into a frame $M$ with complemented image, then there exists a frame map $\bar{f}:\congsL\to M$ such that $f=\bar{f}\circ\nabla$.
\item[(i)] We also note that for frame maps $f$ and $g$ with domain $\congsL$, if $f\circ\nabla=g\circ \nabla$ then $f=g$.

\item[(j)] A useful congruence for our purposes is the \textit{Madden} congruence, denoted $\pi_L$ below: 
\benu
\item[(i)] For $x\in L$, set $P_x=\{t\in L:t\meet x=0\}$. 
\item[(ii)]  For  $x\in L$, $P_x$ is an \sideal, and in $\aichsL$, $P_x=(\darw x)^*$, the pseudocomplement of $\darw x$.
\item[(iii)] Let $\pi_L=\{(x,y):P_x=P_y\}$; $\pi_L$ is an \scong{}. 
\item[(iv)] The quotient map induced by the Madden congruence, $p:L\to L/{\pi_L}$ is dense, onto and the universal such. We refer to this as the \textit{Madden quotient} of $L$. (See \cite{jfasmadden}.)
\eenu
\eenu
\end{definition}

\begin{definition}\label{ab}
For any \sframe{} $L$, define $e_L:\aichsL\to\congsL$ to be the unique frame map such that $e_L(\darw a)=\nabla_a$ for all $a\in L$; that is, making the following diagram commute:
\begin{center}
	\begin{tikzpicture}
		\node at (1.2,2.9){$L$};
		\node at (4.3,2.9){$\aichsL$};
		\node at (1.2,1.3){$\congsL $};
		\draw[<-](1.2,1.6)--(1.2,2.6);
		\draw[->](1.5,2.9)--(3.8,2.9);
		\draw[<-](1.7,1.5)--(3.8,2.6);
		\node at (1,2.1){$\nabla $};
		\node at (2.6,3.2){$\myopdar $};
		\node at (2.8,1.7){$e_L $};
		
	\end{tikzpicture}
\end{center}
\vskip5pt 
That this map $e_L$ exists follows from the freeness of $\aichsL$ as a frame over $L$. See \cite{jfas_hcozforqm}.
\end{definition}

\begin{note}\label{ac}
For any \sframe{} $L$, $\aichsL$ is isomorphic to a subframe of $\congsL$; that is, the free frame over $L$ is isomorphic to a subframe of the frame of \scong{s} of $L$.
\end{note}

\section{Right and left adjoints}

We use the following standard terminology:

\begin{definition}\label{A}
Let $h:L\to M$ be an \sframe{} map.\\
A function $r:M\to L$ is a \textit{right adjoint} of $h$ if
$$h(x)\leq m\dimp x\leq r(m)\textrm{ for all }x\in L, m\in M.$$
A function $l:M\to L$ is a \textit{left adjoint} of $h$ if 
$$l(m)\leq x\dimp m\leq h(x) \textrm{ for all }x\in L, m\in M.$$
\end{definition}

We make no claim that all \sframe{} maps have right (or left) adjoints; this is false (see Example \ref{D}). However, clearly if an \sframe{} map has a right or left adjoint, such is unique.

\begin{lemma}\label{C}
Let $h:L\to M$ be an \sframe{} map.
\benu
\item If $h$ has a right adjoint $r$, then for all $m\in M$,
$$r(m)=\bjoin\{x\in L:h(x)\leq m\}.$$
\item If $h$ has a left adjoint $\ell$, then for all $m\in M$,
$$l(m)=\bmeet\{x\in L:m\leq h(x)\}.$$
\eenu
\end{lemma}

We note that the existence of the above joins and meets has to be established since an \sframe{} need not be complete.

\begin{example}\label{D}
This is an example of an \sframe{} map which has neither a right nor a left adjoint.

Let $L$ be the $\sigma$-frame consisting of all countable and cocountable subsets of $\rnum$, and $\mathbf{2}$ denote the $2$-element chain. Define $h:L\to \mathbf{2}$ by $h(C)=0$ if $C$ is countable and $h(D)=1$ if $D$ is cocountable. Then $h$ is a $\sigma$-frame map. However it has no right adjoint since there is no largest $A\in L$ with $h(A)=0$. Similarly it has no left adjoint. 
\end{example}

\begin{proposition}\label{E}
Let $h:L\to M$ be an \sframe{} map.
\benu
\item Suppose that $h$ has a right adjoint, $r$. Then $h$ preserves all existing joins and $r$ preserves all existing meets.
\item Suppose that $h$ has a left adjoint $L$. Then $h$ preserves all existing meets and $\ell$ preserves all existing joins.
\eenu
\end{proposition}

\section{Closed and open maps}

\begin{definition}\label{I}
Let $h:L\to M$ be an \sframe{} map. 

We call $h$ \emph{closed} if, for all $m \in M$, there exists $x\in L$ with $(h\times h)^{-1}(\nabla_m)=\nabla_x$.

We call $h$ \emph{open} if, for all $m\in M$, there exists $x\in L$ with $(h\times h)^{-1}(\Delta_m)=\Delta_x$.
\end{definition}

We know that (see \cite{jfasmadden}) that $\congs$ is a functor from \sframes{} to frames which is natural in the sense that for any \sframe{} $h:L\to M$ we have a frame map $\congs h:\congsL\to \congsM$ making the following diagram commute:

\begin{center}
	\begin{tikzpicture}
		\node at(1,1.2){$M$}; 	
		\node at(1,3){$L$}; 	
		\node at(4.5,3){$\congsL$};	
		\node at(4.5,1.2){$\congsM$};	
		
		\node at(0.7,2.1){$h$};	
		\node at(2.6,3.3){$\nabla_L $};	
		\node at(2.6,0.8){$ \nabla_M$};	
		\node at(5.2,2.1){$\congs h$};	
		\draw[<-](1,1.5)--(1,2.8); 	
		\draw[->](1.5,3)--(4,3) ;	
		\draw[->](1.5,1.2)--(4,1.2); 	
		\draw[<-](4.7,1.5)--(4.7,2.8);	
	\end{tikzpicture}
\end{center}
\vskip5pt 

Now $(h\times h)^{-1}$ is the right adjoint of $\congs h$, because, for $\theta\in\congsL$, $\congs h(\theta)$ is the \scong{} of $M$ generated by $(h\times h)[\theta]$, so for all $\theta\in\congsL,\phi\in\congsM$,
$$\congs h(\theta)\ctin\phi\dimp\theta\ctin(h\times h)^{-1}(\phi).$$

\begin{theorem}\label{J}
Let $h:L\to M$ be an \sframe{} map.
\begin{enumerate}
	\item[(a)] The map $h$ is closed iff $h$ has a right adjoint, $r$, and for all $x\in L, m\in M$,
	$$r(h(x)\join m)=x\join r(m).$$
	\item[(b)] The map $h$ is open iff $h$ has a left adjoint, $l$, and for all $x\in L,m\in M$,
	$$l(h(x)\meet m)=x\meet l(m).$$
\end{enumerate}
\end{theorem}

\begin{theorem}\label{K}
Let $L$ be an \sframe{} and $\theta$ an \scong{} on $L$.\\
(a) The quotient map $q:L\to L/\theta$ is closed if and only if $\theta$ is a closed \scong; i.e. $\theta=\nabla_a$ for some $a\in L$.\\
(b) The quotient map $q:L\to L/\theta$ is open if and only if $\theta$ is an open \scong{}; that is, $\theta=\Delta_a$ for some $a\in L$.
\end{theorem}

\begin{definition}
Let $h:L\to M$ be an \sframe{} map.  We say that $h$ is \emph{dense} (resp., \emph{codense}) if for all $a\in L, h(a)=0$  (resp., $h(a)=1$) implies that $a=0$ (resp., $a=1$).\end{definition}
\begin{lemma}\label{L} Let $h:L\to M$ be an \sframe{} map.
\benu
\item[(a)] If $h$ is dense and closed, then $h$ is one-one. 
If $h$ is dense, closed and onto, then $h$ is an isomorphism.
\item[(b)] If $h$ is codense and open, then $h$ is one-one.
If $h$ is codense, open and onto, then $h$ is an isomorphism.
\eenu
\end{lemma}

\begin{lemma}\label{M} Suppose that $f:L\to M$ and $g:M\to N$ are \sframe{} maps.
\benu
\item[(a)]
\benu
\item[(i)] If $f$ and $g$ are both closed, then $g\circ f$ is closed. 
\item[(ii)] If $g\circ f$ is closed and $g$ is one-one, then $f$ is closed.
\item[(iii)] If $g\circ f$ is closed and $f$ is onto, then $g$ is closed.
\eenu
\item[(b)] As above but replace ``closed'' by ``open''.
\eenu
\end{lemma}

\section{Boolean properties for partial frames}

The material in this  section comes from \cite{jfasBoolVar}.

We begin by recalling how matters stand in the case of full frames. A Boolean frame is simply a frame that is also a Boolean algebra, that is, every element has a complement. However, Booleanness can also be characterized in a different way. For any frame $M$, let $M_{**}=\{x^{**}:x\in M\}$ where $x^*=\bjoin\{z\in M:z\meet x=0\}$, the pseudocomplement of $x$. The frame map $p:M\to M_{**}$ given by $p(x)=x^{**}$ is called the Booleanization of $M$. It is the least dense quotient of $M$, but is also the unique dense Boolean quotient of $M$. A frame is then Boolean if and only if it is isomorphic to its Booleanization.

Following Madden's lead in \cite{madden}, in \cite{jfasmadden} we constructed a least dense quotient for partial frames. The codomain need not be Boolean, however, as Madden already noted in the case of $\kappa$-frames. We use his terminology, ``d-reduced'', to refer to those partial frames isomorphic to their least dense quotients. We refer the reader to Definition \ref{aaaa}(l) for our notation and terminology in this regard.

The next result characterizes those \sframe{s}, $L$, that are Boolean algebras, in several ways. These involve the free frame over $L$, the congruence frame of $L$ and the the relationship between these two entities.
\begin{proposition}\label{ci}
Let $L$ be an \sframe. The following are equivalent:
\benu
\item $L$ is Boolean; that is, every element of $L$ is complemented.
\item All principal \sideal{s} in $\aichsL$ are complemented.
\item The embedding $e:\aichsL\to\congsL$ is an isomorphism.
\item Every \scong{} $\theta$ of $L$ is an arbitrary join of \scong{s} of the form $\nabla_a$, for some $a\in L$.
\eenu
\end{proposition}

In our experience with partial frames, it has often proved useful to compare properties for a partial frame with the analogous properties for the corresponding free frame. We do this now for Booleanness.

We recall that, if $M$ is a frame and $x\in M$, we call $x$ a \textit{dense} element of $M$ if $x^*=0$.

\begin{proposition}\label{cj}
Let $L$ be an \sframe. The following are equivalent:
\benu
\item The frame $\aichsL$ is Boolean.
\item $\darw 1$ is the only dense element of $\aichsL$.
\item The \sframe{} embedding $\nabla: L\to\congsL$ is an isomorphism.
\item Every $\theta\in\congsL$ has the form $\theta=\nabla_a$, for some $a\in L$.
\eenu
\end{proposition}

We now provide four provably distinct conditions akin to Booleanness for partial frames. In the setting of (full) frames they all amount to every element being complemented.

\begin{theorem}\label{da}
Let $L$ be an \sframe. In the following list of conditions, each one implies the succeeding one, but not conversely.
\benu
\item $\aichsL$ is a Boolean frame.
\item $L$ is a Boolean frame.
\item $L$ is a Boolean \sframe.
\item $L$ is a d-reduced \sframe.
\eenu
\end{theorem}
\bproof
(a)$\imp$(b):
($\not\limp$) See Example \ref{db}.\\
(b)$\imp$(c): 
($\not\limp$): See Example \ref{dc}.\\
(c)$\imp$(d): 
($\not\limp$): See Example \ref{bh}. 
\eproof
\begin{example}\label{bh}
Let $\sels$ designate countable subsets, and consider the $\sigma$-frame $L=\calP_C(\rnum)$, which consists of all countable subsets of $\rnum$ together with $\rnum$ as the top element. Countable join is union, binary meet is intersection.\\
Here $(X,Y)\in \chio$ if and only if, for any countable subset $U$, $U\cap X=\emptyset\dimp U\cap Y=\emptyset$, which makes $X=Y$. So $\chio=\Delta$, which makes $\calP_C(\rnum)$ d-reduced. However, $\calP_C(\rnum)$ is clearly not Boolean.
\end{example}
\begin{example}\label{db}
Let $\sels$ designate countable subsets, and let $\mathcal{L}$ consist of all subsets of $\rnum$. Clearly $\mathcal{L}$ is a Boolean frame. We show that $\aichs\mathcal{L}$ is not Boolean, using condition (d) of Proposition \ref{cj}.\\
Let $\mathcal{I}=\{X\ctin\rnum: X\cap(\rnum\backslash\qnum) \textrm{ is countable}\}$. Then $\mathcal{I}$ is a $\sigma$-ideal of $\mathcal{L}$; that is, a downset closed under countable unions. By Definition \ref{aaaa}(g)(iv), $\bjoinl{X\in\mathcal{I}}\nabla_X=\bigcup\limits_{X\in\mathcal{I}}\nabla_X$ and this cannot have the form $\nabla_Z$ for any $Z\in\mathcal{L}$, since that would require $Z\supseteq X$ for all $X\in\mathcal{I}$, and hence $Z=\rnum$; a contradiction.\phantom{XXXX}
\end{example}

\begin{example}\label{dc}
Let $L$ consist of all countable and co-countable subsets of the real line, and let $\sels$ designate countable subsets. Clearly $L$ is a Boolean $\sigma$-frame, but not a complete lattice, so not a frame.
\end{example}

\section{Comparing congruences on a partial frame and its free frame}

The material in this  section comes from \cite{jfasBoolVar}.

In this section, for a partial frame $L$, we compare $\congsL$, the frame of \scong{s} of $L$, with $\calC(\aichsL)$, the frame of (frame) congruences on $\aichsL$, the free frame over $L$. The universal property of the embedding $\nabla:L\to\congsL$ provides a frame map $E_L:\congsL\to\calC(\aichsL)$. We give an explicit description of this map, and show that it provides a natural transformation.

We then turn our attention to its right adjoint $D_L:\calC(\aichsL)\to \congsL$. Again, we provide an explicit description of this function, including an interesting and useful action on closed congruences (Lemma \ref{ec}).

\begin{definition}\label{dg}
Let $L$ be an \sframe. Consider this diagram
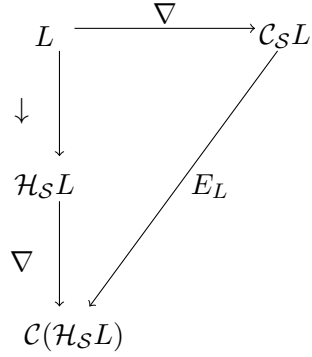
\begin{center}
	\begin{tikzpicture}
		\node at (0.4,1){$\calC(\aichsL)$};
		\node at (0,3){$\aichsL$};
		\node at (0,5){$L$};
		\node at (3.2,5){$\congsL$};
		\node at (-0.3,2){$\nabla$};
		\node at (-0.3,4){$\darw$};
		\node at (1.6,5.3){$\nabla$};
		\node at (2.2,3){$E_L$};
		\draw[->](0.2,4.8)--(0.2,3.4);
		\draw[->](0.2,2.8)--(0.2,1.4);
		\draw[->](0.4,5.1)--(2.8,5.1);
		\draw[->](3.1,4.8)--(0.6,1.4);
	\end{tikzpicture}
\end{center}

By the universal property of $\nabla:L\to\congsL$ there exists a unique frame map $E_L:\congsL\to\calC(\aichsL)$ such that $E_L\circ\nabla=\nabla\circ\darw$; that is, for all $a\in L$ 
$$ E_L(\nabla_a)=\nabla_{\darw a}.$$

\end{definition}

\begin{lemma}\label{dh} Let $L$ be an \sframe{}.
\benu
\item For $\theta\in\congsL$, $E(\theta)$ is the frame congruence on $\aichsL$ generated by $\{(\darw x,\darw y):(x,y)\in\theta\}$; this is denoted by $E(\theta)=\langle (\darw x,\darw y):(x,y)\in\theta\rangle$.
\item The frame map $E:\congsL\to\calC(\aichsL)$ is dense.
\eenu
\end{lemma}

\begin{corollary}\label{di}
Let $L$ be an \sframe{} and $E:\congsL\to\calC(\aichsL)$ given as in Definition \ref{dg}. For all $a\in L$:
\benu
\item $E(\nabla_a)=\nabla_{\darw a}$
\item $E(\Delta_a)=\Delta_{\darw a}$
\eenu
\end{corollary}

\begin{remark}\label{dj}
Let $L$ be an \sframe. the embedding $e:\aichsL\to\congsL$ of Definition \ref{ab} can be incorporated into the diagram of Definition \ref{dg} as follows:
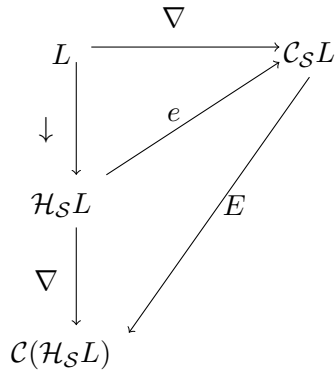
\begin{center}
	\begin{tikzpicture}
		\node at (0,0){$\calC(\aichsL)$};
		\node at (0,2){$\aichsL$};
		\node at (0,4){$L$};
		\draw[->](0.2,3.9)--(0.2,2.4);
		\draw[->](0.2,1.7)--(0.2,0.4);
		\node at (3.3,4){$\congsL$};
		\draw[->](0.4,4.1)--(2.9,4.1);
		\draw[->](0.6,2.4)--(2.9,3.9);
	\draw[->](3.3,3.7)--(0.9,0.3);
		\node at (-0.2,1){$\nabla$};
		\node at (-0.2,3){$\darw$};
		\node at (1.5,4.5){$\nabla$};
		\node at (1.5,3.2){$e$};
		\node at (2.3,2){$E$};
	\end{tikzpicture}
\end{center}
Note that \begin{itemize}
	\item the upper triangle commutes, since $e\circ \darw =\nabla$.
	\item the lower triangle commutes, since, for $I\in\aichsL$,\\ $E\circ e(I)=E(\bjoinl{i\in I}\nabla_i)=\bjoinl{i\in I}E(\nabla_i)=\bjoinl{i\in I}\nabla_{\darw i}=\nabla_I$.\\
	Alternatively, this can be seen because the outer diagram commutes and every \sideal{} is generated by principal \sideal{s}.
\end{itemize}

\end{remark}
\begin{proposition}\label{eg} The function $E_L$ provides a natural transformation from the functor $\calC_{\sels}$ to the functor $\calC\aichs$.
\end{proposition}

We now define, for any \sframe{} $L$, the function $D_L$. In a subsequent lemma, $D_L$ is seen to be the right adjoint of the frame map $E_L$.

\begin{definition}\label{ea} Let $L$ be an \sframe, and $\Phi$ a frame congruence on the frame $\aichsL$. Define $$D_L(\Phi)=\{(x,y)\in L\times L: (\darw x,\darw y)\in\Phi\}.$$

\end{definition}
\begin{lemma}\label{eb}
Let $L$ be an \sframe.
\benu
\item For any frame  congruence $\Phi$ on $\aichsL$, $D_L(\Phi)$ is an \scong{} on $L$.
\item The function $D_L:\calC(\aichsL)\to\congsL$ is the right adjoint of the frame map $E_L:\congsL\to\calC(\aichsL)$ of Definition \ref{dg}.
\item The function $D_L$ preserves bottom, top and arbitrary meets.
\eenu
\end{lemma}

We now provide further properties of $D$, including its action on certain special congruences. We note that the proof of Lemma \ref{ec}(a) uses the fact that, for $I$ an \sideal{} of an \sframe{} $L$, $\bjoinl{i\in I}\nabla_i=\bigcup\limits_{i\in I}\nabla_i$. This is not immediately obvious, but was proved in \cite{jfasheyt} Lemma 3.1.

\begin{lemma}\label{ec}
Let $L$ be an \sframe, and $D$ as in Definition \ref{ea}.
\benu
\item For all $I\in\aichsL$, $D(\nabla_I)=\bcupl{i\in I}\nabla_i.$
\item For all $a\in L$,
\benu
\item $D(\nabla_{\darw a})=\nabla_a$
\item $D(\Delta_{\darw a})=\Delta_a$
\eenu
\item For $I\in\aichsL$, $I$ is principal $\dimp D(\nabla_I)\join D(\Delta_I)=\nabla$.
\eenu
\end{lemma}

\begin{definition} Let $M$ be a full frame. For any $a\in M$ we say $a$ is an \textit{\slind{} element} of $M$ if the following condition holds:

If $a = \bjoin B$ for some $B\ctin M$, then $a =\bjoin D$ for some designated subset $D$ of $M$ such that $D \ctin  B$.
\end{definition}
See \cite{jfas_hcozforqm} for details about this notion. In particular, Lemma 4.3 of that paper characterizes the \slind{} elements of $\aichsL$ as being the principal \sideal{s}.

The next result characterizes those rather special \sframe{s} $L$ for which $E_L$ is an isomorphism.

\begin{theorem}\label{ed}
Let $L$ be an \sframe. The following are equivalent:
\benu
\item The embedding $\darw:L\to\aichsL$ is an isomorphism.
\item Every \sideal{} of $L$ is principal.
\item $L$ is a frame and every element of $L$ is \slind.
\item The frame map $E:\congsL\to\calC(\aichsL)$ is an isomorphism.
\eenu
\end{theorem}

The equivalent conditions of Theorem \ref{ed} might seem rather strong. Here are some examples which show that these can obtain.

\begin{example}\label{ee}
The conditions of Theorem \ref{ed} hold in the following examples:
\begin{itemize}
	\item $\sels$ selects finite subsets and $L$ is a finite frame.
	\item $\sels$ selects countable subsets, and $L$ consists of the open subsets of the real line.
	\item $\sels$ selects finite subsets, or $\sels$ selects countable subsets, and $L$ consists of the cofinite subsets of the real line, together with the empty set.
\end{itemize}
\end{example}

\section{Closed and open embeddings into the free frame and the congruence frame}

\begin{theorem}\label{U}
Let $L$ be an \sframe{} and $\da:L \to \aichsL$ the embedding into its free frame.
\benu
\item[(a)] The map $\da$ has a right adjoint iff $\da$ is an isomorphism.
\item[(b)] The map $\da$ is closed iff $\da$ is an isomorphism.
\item[(c)] The map $\da$ has a left adjoint iff $L$ is a complete lattice.
\item[(d)]  The map $\da$ is open iff $L$ is a frame.
\eenu
\end{theorem}

\begin{theorem}\label{V}
Let $L$ be an \sframe{} and $\nabla:L \to \congsL$ the embedding into its congruence frame.
\benu
\item[(a)] The map $\nabla$ is closed iff $\nabla$ is an isomorphism.
\item[(b)]  The map $\nabla$ is open iff $L$ is a Boolean frame.
\eenu\end{theorem}

\end{document}